# Local inequalities for plurisubharmonic functions

By Alexander Brudnyi*

## Abstract

The main objective of this paper is to prove a new inequality for plurisubharmonic functions estimating their supremum over a ball by their supremum over a measurable subset of the ball. We apply this result to study local properties of polynomial, algebraic and analytic functions. The paper has much in common with an earlier paper [Br] of the author.

## 1. Introduction and formulation of main results

1. A real-valued function $f$ defined on a domain $\Omega \subset \mathbb{C}^n$ is called *plurisubharmonic in* $\Omega$ if $f$ is upper semicontinuous and its restriction to components of a complex line intersected with $\Omega$ is subharmonic.

The main objective of this paper is to prove a new inequality for plurisubharmonic functions estimating their supremum over a ball by supremum over a measurable subset of the ball. The inequality has many applications, several of which are presented in this paper. To formulate the result and its applications we introduce

*Definition* 1.1. A plurisubharmonic function $f : \mathbb{C}^n \longrightarrow \mathbb{R}$ belongs to class $\mathcal{F}_r$ $(r > 1)$ if it satisfies

(i) $$\sup_{B_c(0,r)} f = 0;$$

(ii) $$\sup_{B_c(0,1)} f \geq -1.$$

Hereafter $B(x,\rho)$ and $B_c(x,\rho)$ denote the Euclidean ball with center $x$ and radius $\rho$ in $\mathbb{R}^n$ and $\mathbb{C}^n$, respectively.

*Research supported in part by NSERC.
1991 *Mathematics Subject Classification*. Primary 31B05. Secondary 46E15.

Let the ball $B(x,t)$ satisfy

(1.1) $$B(x,t) \subset B_c(x,at) \subset B_c(0,1),$$

where $a > 1$ is a fixed constant.

THEOREM 1.2. *There are constants $c = c(a,r)$ and $d = d(n)$[1] such that the inequality*

(1.2) $$\sup_{B(x,t)} f \leq c \log\left(\frac{d|B(x,t)|}{|\omega|}\right) + \sup_{\omega} f$$

*holds for every $f \in \mathcal{F}_r$ and every measurable subset $\omega \subset B(x,t)$.*

To illustrate the possible applications of the main result, let us consider a real polynomial $p \in \mathbb{R}[x_1, \ldots, x_n]$ of degree at most $k$ (we will denote the space of these polynomials by $\mathcal{P}_{k,n}(\mathbb{R})$). According to the classical Bernstein "doubling" inequality

(1.3) $$\max_{B_c(0,r)} |p| \leq r^k \max_{B_c(0,1)} |p| \quad (r > 1).$$

Consider the plurisubharmonic function

$$F_r(z) := (\log r)^{-1} k^{-1} (\log |p(z)| - \sup_{B_c(0,r)} \log |p|) \quad (z \in \mathbb{C}^n).$$

From the definition of $\mathcal{F}_r$ and (1.3) it follows that $F_r \in \mathcal{F}_r$ for any $r > 1$. Applying (1.2) with $r = 2$ to this function we get

(1.4) $$\sup_B |p| \leq \left(\frac{d|B|}{|\omega|}\right)^{ck} \sup_{\omega} |p|$$

for an arbitrary ball $B$ and its measurable subset $\omega$.

In fact, in this case, we can take $c = 1$ and $d = 4n$ as follows from the sharp inequality due to Remez [R] for $n = 1$ and Yu. Brudnyĭ-Ganzburg [BG] [2] in the general case.

2. Applications of the main theorem are related to Yu. Brudnyĭ-Ganzburg type inequalities for polynomials, algebraic functions and entire functions of exponential type. We give also applications to log-BMO properties of real analytic functions, which previously were known only for polynomials (see [St]). As is seen from the proof of (1.4) the main result serves in these applications as a kind of amplifier, transforming weak-type inequalities into strong-type ones. Of course, these "weak" inequalities are, clearly, highly nontrivial and obtaining them may require a great deal of effort. Fortunately, a number of these have recently been proved in connection with different aspects of modern analysis (see, in particular, [S], [FN1], [FN2], [FN3], [Br], [BMLT], [RY], [LL]).

---

[1] Here and below the notation $C = C(\alpha, \beta, \gamma, \ldots)$ means that the constant depends only on the parameters $\alpha, \beta, \gamma, \ldots$.

[2] In the original version the ratio on the right-hand side of (1.4) can be replaced by $T_k\left(\frac{1+\beta_n(\lambda)}{1-\beta_n(\lambda)}\right)$ with $\lambda := \frac{|\omega|}{|B|}$ and $\beta_n(\lambda) = (1-\lambda)^{\frac{1}{n}}$. Here $T_k$ is the Tchebychef polynomial of degree $k$.

3. We now formulate two consequences of the main result which give a refinement (and a relatively simple alternative proof) of the basic results of [Br] and [FN3]. We begin with a sharpening of the main result in [Br] (in the original version of inequality (1.5) below the exponent depends on $k$ in a nonlinear way).

To formulate the result suppose that $V \subset \mathbb{R}^n$ is a real algebraic variety of pure dimension $m$ ($1 \leq m \leq n-1$). We endow $V$ with the metric and the measure induced from the Euclidean metric and Lebesgue measure of $\mathbb{R}^n$.

THEOREM 1.3. *For every regular point $x \in V$ there is an open neighborhood $N = N(V)$ of $x$ such that*

$$(1.5) \qquad \sup_B |p| \leq \left(\frac{d\lambda_V(B)}{\lambda_V(\omega)}\right)^{\alpha k \deg(V)} \sup_\omega |p|$$

*for every ball $B \subset N$, measurable subset $\omega \subset B$ and polynomial $p \in \mathcal{P}_{k,n}(\mathbb{R})$.*

*Here $\lambda_V$ denotes the induced Lebesgue measure in $V$ and $d = d(m)$ and $\alpha$ is an absolute constant.*

Our next result is a generalization of the first main result in [FN3] in which $\omega$ in (1.6) below is a ball.

THEOREM 1.4. *Let $F_{1,\lambda}, \ldots, F_{N,\lambda}$ be holomorphic functions on the ball $B_c(0, 1+r_0) \subset \mathbb{C}^n$, $r_0 > 0$, depending real-analytically on $\lambda \in U \subset \mathbb{R}^m$ where $U$ is open. Let $V_\lambda$ be the linear span of the $F_{k,\lambda}, 1 \leq k \leq N$. Then for any compact set $K \subset U$, there is a constant $\gamma = \gamma(K, r_0) > 0$ such that the Yu. Brudnyĭ-Ganzburg type inequality*

$$(1.6) \qquad \sup_{B(x,\rho)} |F| \leq \left(\frac{d(n) |B(x,\rho)|}{|\omega|}\right)^\gamma \sup_\omega |F|$$

*holds for any $F \in V_\lambda, \lambda \in K$ and $\omega \subset B(x, \rho) \subset B(0, 1)$.*

4. Our next results deal with log-BMO properties of algebraic and analytic functions. The estimate of Theorem 1.2 implies BMO-norm estimates for important classes of analytic functions. We formulate only a few results of this kind. Our first result completes Theorem 5.5 of [Br].

THEOREM 1.5. *Let $V$ be a compact algebraic submanifold of $\mathbb{R}^n$. Then for every real polynomial $p \in \mathcal{P}_{k,n}(\mathbb{R})$ with $p|_V \neq 0$ the function $(\log |p|)|_V \in$ BMO$(V)$ and its BMO-norm is bounded above by $C(V)k$.*

Let us recall that the BMO-norm of $f \in L_1(V, d\lambda_V)$ is defined as

$$|f|_* := \sup_B \frac{1}{\lambda_V(B)} \int_B |f - f_B| d\lambda_V,$$

where $f_B := \frac{1}{\lambda_V(B)} \int_B f \, d\lambda_V$, $B \subset V$ is a ball with respect to the induced metric and $\lambda_V$ is the Lebesgue measure on $V$ induced from $\mathbb{R}^n$.

*Remark* 1.6. In the previous version of this result, the BMO-norm was estimated by a constant depending nonlinearly on the degree $k$.

Now let $\{F_{j,\lambda}\}_{1 \leq j \leq N}$ be a family of real analytic functions defined on a compact real analytic manifold $V$ and depending real-analytically on $\lambda$ varying in an open subset $U$ of $\mathbb{R}^m$.

THEOREM 1.7. *Let $V_\lambda := \mathrm{span}\{F_{j,\lambda}\}_{1 \leq j \leq N}$. Then for every compact set $K \subset U$ there is a constant $C = C(K) > 0$ such that*

$$\Big|\log|F|\Big|_* \leq C$$

*for every $F \in V_\lambda$ with $\lambda \in K$.*

## 2. Proof of Theorem 1.2

1. The proof is divided into three parts, the first of which will be presented in this section. It contains several auxiliary results on subharmonic and plurisubharmonic functions.

Let $\mathcal{PSH}(\Omega)$ denote the class of plurisubharmonic in $\Omega$ functions. An important subclass of $\mathcal{PSH}(\mathbb{C}^n)$ is introduced as follows.

*Definition* 2.1. A function $u \in \mathcal{PSH}(\mathbb{C}^n)$ belongs to the class $\mathcal{L}(\mathbb{C}^n)$ (*of functions of minimal growth*) if

(2.1) $$u(z) - \log(1 + |z|) \leq \alpha \quad (z \in \mathbb{C}^n)$$

for a constant $\alpha$.

To formulate our first auxiliary result consider the family $\mathcal{A}_r$ of continuous nonpositive subharmonic functions $f : \mathbb{D} \longrightarrow \mathbb{R}$ such that

(2.2) $$-1 \leq \sup_{\mathbb{D}_r} f.$$

Here $\mathbb{D}_r := \{z \in \mathbb{C}; \ |z| < r\}$, $\mathbb{D} := \mathbb{D}_1$ and $r$ is a fixed number, $0 < r < 1$.

PROPOSITION 2.2. *For every $f \in \mathcal{A}_r$ there exists a subharmonic function $h_f : \mathbb{C} \longrightarrow \mathbb{R}$ and a constant $c_f > 0$ such that*

(i) $$h_f/c_f \in \mathcal{L}(\mathbb{C});$$

(ii) $$f = h_f \text{ on } \mathbb{D}_r;$$

(iii) $$\sup_{f \in \mathcal{A}_r} c_f < \infty.$$

*Proof.* Let $R := \{\frac{1+3r}{4} \leq |z| \leq \frac{1+r}{2}\}$ be an annulus in $\mathbb{D} \setminus \mathbb{D}_r$, and let $X$ denote the family of concentric circles centered at 0 and contained in $R$.

LEMMA 2.3. *Let $f \in \mathcal{A}_r$ and*
$$t(f) := \sup_{S \in X} \inf_{z \in S} f(z).$$

*Then*

(2.3) $$C(r) := \inf_{f \in \mathcal{A}_r} t(f) > -\infty.$$

*Proof.* Below we follow a scheme suggested by N. Levenberg that essentially simplifies our original proof. Let $\{f_i\}_{i \geq 1} \subset \mathcal{A}_r$ be such that
$$\lim_{i \to \infty} t(f_i) = C(r).$$
Without loss of generality we may assume that the sequence does not contain the zero function. For every $S \subset X$ we set
$$S_i := \{z \in S; \ f_i(z) = \min_S f_i\}$$
and
$$K_i := \bigcup_{S \in X} S_i.$$

By the continuity of $f_i$ the set $K_i$ is compact. The set of radii of points in $K_i$ fills out the interval $I_r := [\frac{1+3r}{4}, \frac{1+r}{2}]$ of length $w(r) := \frac{1-r}{4}$. Then the transfinite diameter $\delta(K_i)$ of $K_i$ satisfies

(2.4) $$\delta(K_i) \geq \delta(I_r) = \frac{w(r)}{4}.$$

(See [G, Chap. VII], for the definition and properties of transfinite diameter.)

Now we set

(2.5) $$m_i := \max_{K_i} f_i \quad \text{and} \quad g_i := \frac{f_i}{|m_i|}.$$

Here $m_i < 0$, for otherwise $f_i$ equals 0 identically. To complete the proof we must estimate $|m_i|$ by a constant independent of $i \geq 1$. To this end we will compare $g_i$ with the relative extremal function $u_{K_i, \mathbb{D}}$ of the pair $(K_i, \mathbb{D})$. Recall that the latter is defined by

(2.6) $$u_{K_i, \mathbb{D}}(z) := \sup\{v(z) : \ v \in \mathcal{SH}(\mathbb{D}), v|_{K_i} \leq -1, v \leq 0\}$$

for $z \in \mathbb{D}$. Here $\mathcal{SH}(\mathbb{D}) = \mathcal{PSH}(\mathbb{D})$ for $n = 1$. Since $g_i \leq -1$ on $K_i$ by definition, we have

$$(2.7) \qquad g_i \leq u_{K_i,\mathbb{D}}.$$

Let now $(u_{K_i,\mathbb{D}})^*$ be the upper semicontinuous regularization of (2.6). Then this function is subharmonic in $\mathbb{D}$, see, e.g. [K]. By the nonpositivity of both the regularization and $g_i$ and by inequality (2.7) we have

$$|(u_{K_i,\mathbb{D}})^*| \leq |g_i| = \frac{|f_i|}{|m_i|},$$

as well. From here it follows that at a certain point $z_0 \in \mathbb{D}$, which we will specify later, we get

$$(2.8) \qquad |m_i| \leq \frac{|f_i(z_0)|}{|(u_{K_i,\mathbb{D}})^*(z_0)|}.$$

To select $z_0$ and to estimate the denominator in (2.8) we make use of the relation between the relative extremal function and the capacity $\mathrm{cap}(K_i, \mathbb{D})$ which is defined by

$$\mathrm{cap}(K_i, \mathbb{D}) := \int_{\mathbb{D}} \Delta (u_{K_i,\mathbb{D}})^* dxdy;$$

see, e.g., [K]. Since $(u_{K_i,\mathbb{D}})^*$ satisfies the Laplace equation outside of $K_i$ we can rewrite the right side as follows.

Let $R' \subset \mathbb{D}$ be an arbitrary annulus outside of the circle $\mathrm{conv}(R) = \{z;\ |z| \leq \frac{1+r}{2}\}$ and $\rho$ be a smooth function with support in $\mathrm{conv}(R')$ that equals 1 in $\mathrm{conv}(R') \setminus R'$. Then by Green's formula

$$\mathrm{cap}(K_i, \mathbb{D}) = \int_{\mathbb{D}} \rho \Delta (u_{K_i,\mathbb{D}})^* dxdy = \left| \int_{R'} (u_{K_i,\mathbb{D}})^* \Delta \rho\, dxdy \right| \leq C \max_{R'} |(u_{K_i,\mathbb{D}})^*|.$$

Since the function $(u_{K_i,\mathbb{D}})^*$ is nonpositive and harmonic in $\mathbb{D} \setminus K_i$, Harnack's inequality (see, e.g., [K, Lemma 2.2.9]) implies

$$\max_{R'} |(u_{K_i,\mathbb{D}})^*| \leq C' |(u_{K_i,\mathbb{D}})^*(z_0)|$$

for a constant $C'$ (depending on $r$ only) and every $z_0 \in R'$.

Putting together (2.8) and the latter two inequalities we find that the inequality

$$|m_i| \leq \frac{C''|f_i(z_0)|}{\mathrm{cap}(K_i, \mathbb{D})}$$

holds for every $z_0 \in R'$. But by the definition of $\mathcal{A}_r$ and the maximum principle, $0 > \max_{R'} f_i \geq -1$. Taking $z_0$ as a point at which the latter maximum is attained, we then get

$$|m_i| \leq \frac{C''}{\mathrm{cap}(K_i, \mathbb{D})}.$$

It remains to apply the one-dimensional version of the *comparison theorem* of Alexander and Taylor, see [AT], that gives the following inequality for the transfinite diameter of $K_i$ (which coincides with the logarithmic capacity of $K_i$):

$$\delta(K_i) \leq \exp\left(-\frac{2\pi}{\operatorname{cap}(K_i, \mathbb{D})}\right).$$

Putting together the latter two inequalities and inequality (2.4) we finally obtain

$$|m_i| \leq C''' \log\left(\frac{4}{w(r)}\right) =: C'''(r)$$

for every $i \geq 1$. By the definition of the $\{f_i\}$ it follows that

$$\inf_{f \in \mathcal{A}_r} t(f) = \lim_{i \to \infty} t(f_i) \geq -\inf_i |m_i| \geq -C'''(r) > -\infty.$$

The lemma is proved. □

Now we are in a position to prove Proposition 2.2. Let $f \in \mathcal{A}_r$. According to the lemma there exists a circle $S_f \in X$ such that

$$\inf_{S_f} f \geq C(r) > -\infty.$$

$S_f$ is the boundary of the disk $\mathbb{D}_{r(f)}$, where

$$\frac{1+3r}{4} \leq r(f) \leq \frac{1+r}{2}.$$

We now define the required subharmonic function $h_f(z) : \mathbb{C} \longrightarrow \mathbb{R}$ by

$$h_f(z) := \begin{cases} f(z) & (z \in \mathbb{D}_{r(f)}) \\ \max\left\{f(z), \dfrac{2C(r) \log \frac{4|z|}{3+r}}{\log \frac{4r(f)}{3+r}}\right\} & (z \in \mathbb{D} \setminus \mathbb{D}_{r(f)}) \\ \dfrac{2C(r) \log \frac{4|z|}{3+r}}{\log \frac{4r(f)}{3+r}} & (z \in \mathbb{C} \setminus \mathbb{D}). \end{cases}$$

Since the ratio in the third formula is less than $C(r) < 0$ on $S_f$ and greater than 0 on $\partial \mathbb{D}$, and since $f$ is continuous, $h_f$ is subharmonic on $\mathbb{C}$. Moreover, according to Definition 2.1,

$$\frac{\log \frac{4r(f)}{3+r}}{2C(r)} h_f \in \mathcal{L}(\mathbb{C}).$$

It remains to define

$$c_f := \frac{\log \frac{4r(f)}{3+r}}{2C(r)}.$$

Then

$$c_f \leq \frac{\log \frac{1+3r}{3+r}}{2C(r)} < \infty$$

and the proposition is proved. □

The final result of this section discusses an approximation theorem for plurisubharmonic functions which will allow us to reduce the proof to the case of $C^\infty$-functions.

Let $\kappa$ be a nonnegative radial $C^\infty$-function on $\mathbb{C}^n$ satisfying

$$(2.9) \qquad \int_{\mathbb{C}^n} \kappa(x)dxdy = 1, \quad \mathrm{supp}(\kappa) \subset B_c(0,1),$$

where $z = x + iy$, $x, y \in \mathbb{R}^n$. Let $\Omega \subset \mathbb{C}^n$ be a domain. For $f \in \mathcal{PSH}(\Omega)$, we let $f_\varepsilon$ denote the function defined by

$$(2.10) \qquad f_\varepsilon(w) := \int_{\mathbb{C}^n} \kappa(z) f(w - \varepsilon z) dx dy,$$

where $w \in \Omega_\varepsilon := \{z \in \Omega : \mathrm{dist}(z, \partial\Omega) > \varepsilon\}$. It is well known, see, e.g., [K, Th. 2.9.2], that $f_\varepsilon \in C^\infty \cap \mathcal{PSH}(\Omega_\varepsilon)$ and that $f_\varepsilon(w)$ monotonically decreases and tends to $f(w)$ for each $w \in \Omega$ as $\varepsilon \to 0$.

LEMMA 2.4. *Let $f \in \mathcal{F}_r$. Assume that the functions $\{f_{1/k}\}_{k \geq k_0}$ satisfy inequality (1.2) with $B(x,t)$ and a compact $\omega$ independent of $k$. Then $f$ also satisfies this inequality.*

*Proof.* Since $f$ is defined on $B_c(0,r)$ with $r > 1$, the function $f_{1/k}$ belongs to $C^\infty \cap \mathcal{PSH}(B_c(0, r - 1/k))$. Let $\{w_k\}_{k \geq 1} \subset \omega$ be such that

$$f_{1/k}(w_k) = \max_\omega f_{1/k}.$$

Moreover, let $z_{\varepsilon,t} \in B(x,t)$ be a point such that

$$(2.11) \qquad \sup_{B(x,t)} f - f(z_{\varepsilon,t}) < \varepsilon.$$

According to the assumptions of the lemma

$$(2.12) \quad f(z_{\varepsilon,t}) = \lim_{k \to \infty} f_{1/k}(z_{\varepsilon,t}) \leq c \log \frac{d|B(x,t)|}{|\omega|} + \limsup_{k \to \infty} f_{1/k}(w_k).$$

To estimate the second summand let us use (2.9) and (2.10):

$$(2.13) \qquad f_{1/k}(w_k) = \int_{\mathbb{C}^n} \kappa(z) f(w_k - z/k) dx dy \leq \sup_{B_c(w_k, 1/k)} f.$$

Now let $x_k \in B_c(w_k, 1/k)$ be such that the supremum on the right is less than $f(x_k) + 1/k$. Because of the compactness of $\omega$ we can assume that $w := \lim_{k \to \infty} w_k$ exists. Then we have $\lim_{k \to \infty} x_k = \lim_{k \to \infty} w_k = w \in \omega$. Using the upper semicontinuity of $f$ it follows that

$$\limsup_{k \to \infty} \sup_{B_c(w_k, 1/k)} f \leq \limsup_{k \to \infty} f(x_k) \leq f(w) \leq \sup_\omega f,$$

which leads to the inequality

$$\limsup_{k \to \infty} f_{1/k}(w_k) \leq \sup_\omega f.$$

Putting this inequality together with (2.11) and (2.12) and letting $\varepsilon \to 0$, we get
$$\sup_{B(x,t)} f \leq c \log \frac{d|B(x,t)|}{|\omega|} + \sup_{\omega} f.$$
The proof is complete. □

2. The second part of the proof of Theorem 1.2 is the Bernstein "doubling" inequality for functions in $\mathcal{F}_r$.

PROPOSITION 2.5. *Let $f \in \mathcal{F}_r$ and $s \in [1, a]$, $a > 1$. Suppose that*
$$B_c(x, t) \subset B_c(x, at) \subset B_c(0, 1). \tag{2.14}$$
*Then there is a constant $c = c(r)$ such that*
$$\sup_{B_c(x, st)} f \leq c \log s + \sup_{B_c(x, t)} f. \tag{2.15}$$

*Proof.* Consider the pair of embedded balls $B_c(x, \frac{r-|x|}{r}) \subset B_c(x, r-|x|)$, where $x \in B_c(0, 1)$ and $|\cdot|$ denotes the Euclidean norm. The smaller ball contains $B_c(x, 1-|x|)$ which has maximal radius of balls in $B_c(x, 1)$ centered at $x$. From (2.14) it follows that
$$B_c(x, at) \subset B_c(x, \frac{r-|x|}{r}). \tag{2.16}$$
Let
$$\gamma_r(f; x) := \sup_{B_c(x, \frac{r-|x|}{r})} f$$
and
$$\gamma_r := \inf_{f \in \mathcal{F}_r} \inf_{x \in B_c(0,1)} \gamma_r(f; x). \tag{2.17}$$
Clearly, if $r_1 \leq r_2$ then
$$\gamma_{r_1}(f; x) \leq \gamma_{r_2}(f; x) \quad \text{and} \quad \gamma_{r_1} \leq \gamma_{r_2}. \tag{2.18}$$

LEMMA 2.6. *There is a nonpositive constant $C = C(r)$ such that the inequality*
$$\gamma_r \geq \lim_{k \to \infty} \gamma_{r_k} \geq C > -\infty$$
*holds for every $\{r_k\}_{k \geq 1}$ increasing to $r$.*

*Proof.* According to (2.18), $\{\gamma_{r_k}\}_{k \geq 1}$ is a monotone nondecreasing sequence and therefore $\lim_{k \to \infty} \gamma_{r_k} (\in [-\infty, 0])$ does exist. Let $\{f_k \in \mathcal{F}_{r_k}\}_{k \geq 1}$ and $\{x_k\}_{k \geq 1} \subset B_c(0, 1)$ be chosen such that
$$\lim_{k \to \infty} \gamma_{r_k}(f_k; x_k) = \lim_{k \to \infty} \gamma_{r_k}.$$

Here we may assume that each $f_k$ does not identically equal 0. Let $B_k$ denote the ball $B_c(x_k, \frac{r_k - |x_k|}{r_k})$ and $\lambda B$ denote the homothety of $B$ with center 0 and dilation coefficient $\lambda > 0$. Consider the sequence of balls

$$\{t_k B_k\}_{k \geq 1}, \quad \text{where} \quad t_k := (r/r_k) > 1$$

and the sequence of functions

$$f'_k(z) := f_k(z/t_k) \quad (z \in B_c(0, r)).$$

Then we have

$$\gamma_{r_k}(f_k; x_k) = \sup_{t_k B_k} f'_k.$$

Without loss of generality we assume that $\{t_k x_k\}_{k \geq 1}$ converges to $x \in B_c(0,1)$. Then the limit ball $B_c(x, \frac{r-|x|}{r})$ of the sequence $\{t_k B_k\}$ has radius at least $l := \frac{r-1}{r}$. Therefore its intersection with $B_c(0,1)$ contains the ball $B^l := B_c(y, l/4)$, where $y = x(1 - \frac{r-|x|}{2r|x|})$. Passing to a subsequence we may assume that

$$B^l \subset t_k B_k$$

for all $k \geq 1$. It follows that

(2.19) $$m_k := \sup_{B^l} f'_k \leq \gamma_{r_k}(f_k; x_k) < 0.$$

Consider now the sequence $\{f''_k := f'_k/|m_k|\}_{k \geq 1}$. Each function of the sequence is less than or equal to $-1$ on $B^l$ and nonpositive on $B_c(0, r)$. Therefore it is bounded above by the relative extremal function

(2.20) $$u_{B^l, B_c(0,r)} := \sup\{v(z) : v \in \mathcal{PSH}(B_c(0,r)), v|_{B^l} \leq -1, v \leq 0\}.$$

Since the compact ball $B^l$ is *pluriregular*, this function is continuous and strictly negative outside $B^l$ (see, e.g., [K, Cor. 4.5.9]). Therefore

(2.21) $$M(r_k) := \max_{\partial B_c(0, t_k)} u_{B^l, B_c(0, r)} < 0$$

and

$$|f''_k(z)| \geq |M(r_k)|$$

for every $z \in \partial B_c(0, t_k)$. From this, the definition of $f''_k$ and inequality (2.19), it follows that

$$|f'_k(z)| \geq |M(r_k) \gamma_{r_k}(f_k; x_k)| \quad (z \in \partial B_c(0, t_k)).$$

But the supremum of $f'_k$ over $B_c(0, t_k)$ is at least $-1$. So the previous inequality yields

$$|\gamma_{r_k}(f_k; x_k)| \leq \frac{1}{|M(r_k)|}.$$

Letting $k \to \infty$ we conclude that

(2.22) $$|\lim_{k\to\infty} \gamma_{r_k}(f_k; x_k)| \leq \frac{1}{|M(r)|},$$

where

(2.23) $$M(r) := \max_{\partial B_c(0,1)} u_{B^l, B_c(0,r)} < 0.$$

The proof of the lemma is complete. □

We now proceed to prove Proposition 2.5. Let $\{f_{1/k}\}_{k\geq 1}$ be the approximating sequence of Lemma 2.4 generated by $f \in \mathcal{F}_r$. Since $f \leq 0$ in $B_c(0,r)$ and the smoothing kernel $\kappa$ is a nonnegative function, $f_{1/k} \leq 0$ in $B_c(0, r-1/k)$. Moreover, $\{f_{1/k}(z)\}_{k\geq 1}$ converges monotonically to $f(z)$ at any $z \in B_c(0,r)$. Then for $k$ sufficiently large, say $k \geq k_0$,

$$\sup_{B_c(0,1)} f_{1/k} \geq \sup_{B_c(0,1)} f \geq -1.$$

Thus $f_{1/k} \in C^\infty \cap \mathcal{F}_{r-1/k}$, $k \geq k_0$. We now set $r_k := r - 1/k$ and consider the sequence

$$\left\{ \sup_{B_c(x, \frac{r_k - |x|}{r_k})} f_{1/k} \right\}_{k \geq k_0}.$$

From Lemma 2.6 it follows that

$$\sup_{B_c(x, \frac{r_k - |x|}{r_k})} f_{1/k} > 2C(r)$$

for $k \geq k_1 (\geq k_0)$. Since $f_{1/k} \in C^\infty \cap \mathcal{PSH}(B_c(0, r_k))$, there is a point $z \in \partial B_c(x, \frac{r_k - |x|}{r_k})$ where the supremum is attained. Further, there is an open neighborhood $U$ of $z$ where $f_{1/k}$ is greater than $2C(r)$. Thus there exists a finite family $\mathcal{G}$ of rotations (unitary transformations) of $\mathbb{C}^n$ centered at $x$ such that $\{g(U)\}_{g \in \mathcal{G}}$ forms a covering of an open neighborhood $W$ of $\partial B_c(x, \frac{r_k - |x|}{r_k})$. The plurisubharmonic function

$$g_k(x) := \max_{g \in \mathcal{G}} f_{1/k}(gx) \quad (k \geq k_1)$$

satisfies

(i) $$\max_{B_c(x,s)} g_k = \max_{B_c(x,s)} f_{1/k}$$

for any $s \in (0, r_k - |x|)$ and

(ii) $$g_k(w) > 2C(r)$$

for any $w \in \partial B_c(x, \frac{r_k - |x|}{r_k})$.

Now set
$$c_k := \frac{2C(r)}{\log \frac{r+1}{2r_k}}.$$

Since $r > 1$, this constant is greater than 0 for $k$ sufficiently large, and we can assume that it holds for $k \geq k_1$. We now define for $k \geq k_1$ the function

$$h_k(z) := \begin{cases} g_k(z) & (z \in B_c(x, \frac{r_k - |x|}{r_k})) \\ \max\{c_k \log \frac{(r+1)|z-x|}{2(r_k - |x|)}, g_k(z)\} & (z \in B_c(x, r_k - |x|) \setminus B_c(x, \frac{r_k - |x|}{r_k})) \\ c_k \log \frac{(r+1)|z-x|}{2(r_k - |x|)} & (z \notin B_c(x, r_k - |x|)). \end{cases}$$

Since the logarithmic term is less than $g_k$ on $\partial B_c(x, \frac{r_k - |x|}{r_k})$, by (ii) above, and more than 0 on $\partial B_c(x, r_k - |x|)$ since $r > 1$, the function $h_k$ is plurisubharmonic on $\mathbb{C}^n$. Furthermore, the function $h'_k = \frac{1}{c_k} h_k$ clearly belongs to $\mathcal{L}(\mathbb{C}^n)$. Compare now $h'_k$ with the global $\mathcal{L}$-extremal function $E_{B_c(x,t)}$ defined by

(2.24) $\quad E_{B_c(x,t)}(z) := \sup\{u(z) : u \in \mathcal{L}(\mathbb{C}^n), u \leq 0 \text{ on } B_c(x,t)\}.$

From this definition it follows that

(2.25) $\quad h'_k - \sup_{B_c(x,t)} h'_k \leq E_{B_c(x,t)}.$

Now we make use of the important representation of $\mathcal{L}$-extremal functions (see, e.g., [K, Th. 5.1.7]),

(2.26) $\quad E_{B_c(x,t)}(z) := \sup \left\{ \frac{\log |P(z)|}{\deg P}; \max_{B_c(x,t)} |P| \leq 1 \right\},$

where the supremum is taken over all holomorphic polynomials on $\mathbb{C}^n$. This representation and the classical Bernstein inequality on the growth of univariate complex polynomials lead to the estimate

$$\sup_{B_c(x,st)} E_{B_c(x,t)} \leq \log s \quad (1 \leq s).$$

From (2.25) it follows that

$$\sup_{B_c(x,st)} h_k - \sup_{B_c(x,t)} h_k \leq c_k \log s.$$

From property (i) of $g_k$ and the definition of $h_k$, the doubling inequality (2.15) is obtained for $f_{1/k}$ with the constant $c_k$, $k \geq k_1$. Applying Lemma 2.4 we obtain then the doubling inequality for $f \in \mathcal{F}_r$ with $c = \frac{2C(r)}{\log \frac{r+1}{2r}}$. $\square$

3. *Proof of the Yu. Brudnyĭ-Ganzburg type inequality.* We have to prove that if $f \in \mathcal{F}_r$ and $\omega$ is a measurable subset of $B(x,t) (\subset B_c(x,at) \subset B_c(0,1))$ of positive measure then

$$\text{(2.27)} \qquad \sup_{B(x,t)} f \leq c \log \frac{d|B(x,t)|}{|\omega|} + \sup_\omega f.$$

It is clear that we may assume that $\omega$ is compact. In fact, otherwise $\omega = \omega_0 \cup \left(\bigcup_{j=1}^\infty \omega_j\right)$, where $|\omega_0| = 0$ and $\{\omega_j\}$ is an increasing sequence of compact sets. If (2.27) holds for every $\omega_j$, then we obtain the result for $\omega$ by letting $j \to \infty$ since $c = c(a,r)$ and $d = d(n)$ do not depend on $\omega$.

Taking into account Proposition 2.5 and the approximation of Lemma 2.4, it suffices to prove the following equivalent statement. Let $B_c(0,1) \subset B_c(0,a)$ and $\mathcal{R}_a$ be the family of continuous plurisubharmonic functions $f: B_c(0,a) \longrightarrow \mathbb{R}$ satisfying

(i) $$\sup_{B_c(0,a)} f = 0$$

(ii) $$\sup_{B_c(0,1)} f \geq -c \log a,$$

with the constant $c$ from Proposition 2.5. Then for every measurable subset $\omega \subset B(0,1)$ of positive measure and every $f \in \mathcal{R}_a$,

$$\text{(2.28)} \qquad \sup_{B(0,1)} f \leq c' \log \frac{d|B(0,1)|}{|\omega|} + \sup_\omega f.$$

Here $d = d(n)$ and $c' = c'(a,c)$. In fact, by a translation and a dilation with coefficient $\frac{1}{t}$ we can transform the balls $B_c(x,t)$ and $B_c(x,at)$ into the balls $B_c(0,1)$ and $B_c(0,a)$, respectively. Then the first term on the right in (2.27) does not change. Moreover, the inequality of Proposition 2.5 states that the pullback of a function $f \in \mathcal{F}_r$ determined by this transformation will satisfy to conditions (i) and (ii). Finally, according to Lemma 2.4 we can assume that $f$ is continuous on $B_c(0,1)$.

It remains to prove (2.28). We begin with

LEMMA 2.7.   *There is a constant $C = C(c,a) > 0$ such that*

$$\text{(2.29)} \qquad \max_{B(0,1)} f \geq -C$$

*for every $f \in \mathcal{R}_a$.*

*Proof.* We can repeat the arguments of Lemma 2.6 related to the use of the relative extremal function (2.20). In this case the ball $B(0,1)$ is nonpluripolar. So using the inequality from (ii) we obtain (2.29) with, e.g., $C = \frac{c \log a}{|M(a)|}$, where

$$M(a) := \sup_{\partial B_c(0,(1+a)/2)} u_{B(0,1), B_c(0,a)}.$$

□

Now let $f \in \mathcal{R}_a$ and $x_f \in B(0,1)$ be such that

(2.30) $$M_f := f(x_f) = \max_{B(0,1)} f.$$

By Lemma 3 of [BG] there is a ray $l_f$ with origin at $x_f$ such that

(2.31) $$\frac{\mathrm{mes}_1(B(0,1) \cap l_f)}{\mathrm{mes}_1(\omega \cap l_f)} \leq \frac{n|B(0,1)|}{|\omega|}.$$

Let $l'_f$ be the one-dimensional affine complex line containing $l_f$ and let $z_f$ be a point of $B_c(0,1) \cap l'_f$ such that

$$t_f := \mathrm{dist}(0, l'_f) = |z_f|.$$

Consider the disks

$$\widetilde{D}_f := \frac{1}{r_f}(D_f - z_f) \quad \text{and} \quad \widetilde{D}'_f := \frac{1}{r_f}(D'_f - z_f),$$

where we set

$$D_f := l'_f \cap B_c(0,a), \quad D'_f := l'_f \cap B_c(0,(a+1)/2).$$

The latter sets are disks of radii

$$r_f := \sqrt{a^2 - t_f^2}, \quad r'_f := \sqrt{\left(\frac{a+1}{2}\right)^2 - t_f^2},$$

respectively centered at $z_f$. Note also that $x_f \in D'_f$. Without loss of generality we can identify $\widetilde{D}'_f \subset \widetilde{D}_f$ with the pair $\mathbb{D}_{s_f} \subset \mathbb{D}_1$, where $\mathbb{D}_r := \{z \in \mathbb{C}; |z| \leq r\}$ and $s_f := r'_f/r_f$. The pullback of the restriction $f|_{D_f}$ to $\mathbb{D}_1$ is denoted by $f'$.

LEMMA 2.8. *There exists a number $r = r(a) < 1$ such that*

$$\mathbb{D}_{s_f} \subset \mathbb{D}_r \subset \mathbb{D}_1$$

*for any $f \in \mathcal{R}_a$.*

*Proof.* It follows from the inequality

$$\frac{r'_f}{r_f} = \sqrt{\frac{(a+1)^2 - 4t_f^2}{4(a^2 - t_f^2)}} \leq \frac{a+1}{2a} < 1$$

that one can choose $r(a) = \frac{a+1}{2a}$.

□

In what follows it is worth noting that

$$\max_{\mathbb{D}_r} f' \geq \max_{\mathbb{D}_{s_f}} f' \geq M_f \geq -C$$

by (2.30) and (2.29). We apply now Proposition 2.2 to the function $\frac{f'}{C}$ which clearly satisfies condition (2.2), i.e., belongs to $\mathcal{A}_r$. Returning to the function $f$ we can formulate the result of this proposition as follows.

*Statement.* There exists a subharmonic function $h_f$ defined on $l'_f$ and a constant $c_f > 0$ such that $h_f/c_f$ is of minimal growth on $l'_f$ and $f = h_f$ on $D'_f$. Moreover,

$$c' = c'(c, a) := \sup_{f \in \mathcal{R}_a} c_f < \infty.$$

Set $m_f(\omega) := \max_{\omega \cap l_f} f$ and consider the function $\frac{h_f - m_f(\omega)}{c_f}$. This function is clearly less than or equal to the $\mathcal{L}$-extremal function $E_{\omega \cap l_f}$ (see (2.24) for the definition, replacing $\mathcal{L}(\mathbb{C}^n)$ by $\mathcal{L}(l_f)$ and $B_c(x, t)$ by $\omega \cap l_f$). Using the polynomial representation (2.26) of $\mathcal{L}$-extremal functions and the one-dimensional Yu. Brudnyĭ-Ganzburg inequality (see the remark after (1.4) and (2.31)), we obtain

$$\max_{B(0,1) \cap l_f} E_{\omega \cap l_f} \leq \log \frac{4\mathrm{mes}_1(B(0,1) \cap l_f)}{\mathrm{mes}_1(\omega \cap l_f)} \leq \log \frac{4n|B(0,1)|}{|\omega|}.$$

From the selection of $l_f$ and the statement above, we obtain

$$\max_{B(0,1)} f - \max_\omega f \leq \max_{B(0,1) \cap l_f} h_f - m_f(\omega) \leq c_f \max_{B(0,1) \cap l_f} E_{\omega \cap l_f} \leq c' \log \frac{4n|B(0,1)|}{|\omega|}.$$

This proves the desired inequality (2.28) for $f \in \mathcal{R}_a$ and, hence, Theorem 1.2. □

## 3. Proof of consequences

1. *Proof of Theorem* 1.3. Let $V \subset \mathbb{R}^n$ be a real algebraic variety of pure dimension $m$, $1 \leq m \leq n-1$. Let $x \in V$ be a regular point of $V$. We have to prove that there is an open neighbourhood $N \subset V$ of $x$ depending on $V$ such that

(3.1) $$\sup_B |p| \leq \left(\frac{d\lambda_V(B)}{\lambda_V(\omega)}\right)^{\alpha k \deg(V)} \sup_\omega |p|$$

for every ball $B \subset N$, measurable subset $\omega \subset B$ and real polynomial $p$ of degree $k$. Here $\deg(V)$ is the degree of the manifold $V$, $d = d(m)$ depends only on $\dim V$ and $\alpha$ is an absolute constant.

For the proof we need an estimate for a $p$-valent function regular in $\mathbb{D}_R := \{z \in \mathbb{C};\ |z| < R\}$. This estimate is due to Roytwarf and Yomdin (see [RY, Th. 2.1.3 and Cor. 2.3.1]). We give, for the sake of completeness, a relatively simple proof of the result.

LEMMA 3.1.  *Suppose that $f : \mathbb{D}_R \longrightarrow \mathbb{C}$ is regular and assumes no value more than $p$ times. Then for any $R' < R$ and any $\alpha \in (0,1)$ the inequality*

$$(3.2) \qquad \max_{\mathbb{D}_{R'}} |f| \leq C^p \max_{\mathbb{D}_{\alpha R'}} |f|$$

*holds with $C = C(\alpha, \beta)$, where $\beta := \frac{R'}{R}$.*

*Proof.* We set

$$M := \max_{\mathbb{D}_{R'}} |f| \quad \text{and} \quad M_\alpha := \max_{\mathbb{D}_{\alpha R'}} |f|.$$

If $f(z) = \sum_{k=0}^{\infty} a_k z^k$ then

$$M \leq \sum_{k=0}^{p} |a_k|(R')^k + \sum_{k=p+1}^{\infty} |a_k|(R')^k =: \sum_1 + \sum_2.$$

From the Cauchy inequality for $\mathbb{D}_{\alpha R'}$ we get

$$(3.3) \qquad \sum_1 \leq M_\alpha \sum_{k=0}^{p} \left(\frac{R'}{\alpha R'}\right)^k \leq \frac{\alpha^{-p}}{1-\alpha} M_\alpha.$$

To estimate the second sum we apply the coefficient inequality for $p$-valent functions (see [H]). According to it,

$$(3.4) \qquad |a_j| R^j \leq \left(\frac{A}{p}\right)^{2p} j^{2p} \max_{1 \leq j \leq p} |a_j| R^j$$

for every $j > p$, where $A$ is an absolute constant.

By the Cauchy inequality the maximum on the right of (3.4) can be estimated by

$$M_\alpha \max_{1 \leq j \leq p} \frac{R^j}{(\alpha R')^j} \leq M_\alpha (\alpha \beta)^{-p}.$$

Putting this and (3.4) together we obtain

$$\sum_2 \leq \left(\frac{A}{p}\right)^{2p} (\alpha \beta)^{-p} M_\alpha \sum_{j=p+1}^{\infty} j^{2p} \left(\frac{R'}{R}\right)^j \leq \left(\frac{A}{p}\right)^{2p} (\alpha \beta)^{-p} \phi_{2p}(\beta) M_\alpha,$$

where

$$\phi_l(\beta) = \sum_{k=1}^{\infty} k^l \beta^k.$$

We will prove later that

$$\phi_{2p}(\beta) < \frac{(4p)^{2p}}{(1-\beta)^{2p+1}}, \tag{3.5}$$

which together with the previous inequality leads to the estimate

$$\sum_2 \leq (4A)^{2p}(\alpha\beta)^{-p}(1-\beta)^{-2p-1}M_\alpha.$$

From this and (3.3), the required inequality (3.2) follows with

$$C(\alpha,\beta) = 2\max\left\{\frac{1}{\alpha(1-\alpha)}, \frac{(4A)^2}{\alpha\beta(1-\beta)^3}\right\}.$$

It remains to prove (3.5). To this end, notice that

$$\phi_l(\beta) = \left(\beta\frac{d}{d\beta}\right)^l\left(\frac{1}{1-\beta}\right).$$

Then by induction on $l$ we have

$$\phi_l(\beta) = \frac{p_l(\beta)}{(1-\beta)^{l+1}}, \tag{3.6}$$

where $p_l$ is a polynomial of degree $l$. Moreover, we have the identity

$$p_{l+1}(\beta) = \beta(1-\beta)p'_l(\beta) + (l+1)p_l(\beta).$$

Let $\mu(l) := \max_{0\leq\beta\leq 1}|p_l(\beta)|$. Then from the previous identity and the Bernstein polynomial inequality we get

$$\mu(l+1) \leq \max_{0\leq\beta\leq 1}|\beta(1-\beta)p'_l(\beta)| + (l+1)\mu(l) \leq l\mu(l) + (l+1)\mu(l) = (2l+1)\mu(l).$$

This recurrence inequality implies that

$$\mu(2p) \leq \mu(0)\prod_{l=0}^{2p-1}(2l+1) = \prod_{l=0}^{2p-1}(2l+1) < (4p)^{2p},$$

which combined with (3.6) gives the required estimate (3.5). □

We now proceed to the proof of Theorem 1.3. Let $V_c$ denote the complexification of $V$, i.e., the minimal complex algebraic subvariety of $\mathbb{C}^n$ such that $V$ is a connected component of $V_c \cap \mathbb{R}^n$. Then the regularity of $x$ in $V$ implies that it is a regular point of $V_c$. We will assume without loss of generality that $x$ coincides with the origin $0 \in \mathbb{R}^n$. Then there exist open neighborhoods $U_x \subset U'_x \subset V_c$ of the point $x$ and a linear holomorphic projection of $\mathbb{C}^n$ onto $\mathbb{C}^m$ whose restriction $\phi : V_c \longrightarrow \mathbb{C}^m \subset \mathbb{C}^n$ is biholomorphic in a neighborhood of $U'_x$ such that

(i)   $\phi(x) = 0$,  $\phi(U_x) = B_c(0,r)$, and $\phi(U'_x) = B_c(0,2r)$ for some $r > 0$;
(ii)  $\phi|_V : V \longrightarrow \mathbb{R}^m$

(see, e.g., [Br, Prop. 2.8]). According to (ii), $\phi$ has a smooth inverse defined on $B(0, r)$ and $N' := \phi^{-1}(B(0, r))$ is an open neighborhood of $x$ in $V$.

Consider now a real polynomial $p \in \mathcal{P}_{k,n}$ and its extension $p_c$ to $\mathbb{C}^n$ as a holomorphic polynomial of degree $k$. Let $z \in B_c(0, r)$ be a maximum point of $|p_c \circ \phi^{-1}|$ in $B_c(0, r)$. Let $L$ be a complex line passing through $z$ and the origin. The restriction of $p_c \circ \phi^{-1}$ to $L$ is an algebraic function of one variable, whose local valency can be estimated by the multidimensional Bezout theorem. By this theorem the function $p_c \circ \phi^{-1}$ assumes no value more than $p := k \deg(V_c)$ times in $B_c(0, 2r) \cap L$. Applying Lemma 3.1 to this function we obtain

$$(3.7) \qquad M_r := \max_{B_c(0, 3r/2)} |p_c \circ \phi^{-1}| \leq C^{k \deg(V)} \max_{B_c(0, r)} |p_c \circ \phi^{-1}|,$$

where $C$ is an absolute constant.

Consider now the function $f : B_c(0, 3/2) \longrightarrow \mathbb{R}$ defined by

$$f(z) := \frac{1}{k \log c \deg(V)} (\log |(p_c \circ \phi^{-1})(rz)| - \log M_r).$$

Then $\sup_{B_c(0, 3/2)} f = 0$ and, by (3.7), $\sup_{B_c(0, 1)} f \geq -1$. So $f \in \mathcal{F}_{3/2}$ and the conditions of Theorem 1.2 are fulfilled. Applying Theorem 1.2 to $f$ we get

$$(3.8) \qquad \sup_{\phi(B)} |p \circ \phi^{-1}| \leq \left( \frac{d(m)|\phi(B)|}{|\phi(\omega)|} \right)^{\alpha k \deg(V)} \sup_{\phi(\omega)} |p \circ \phi^{-1}|$$

for every ball $\phi(B) \subset B(0, r)$ and measurable subset $\omega \subset B$. Here $\alpha := c \log C$, where $c$ and $d(m)$ are the constants in the inequality of Theorem 1.2 and $C$ is that of (3.7). To finish the proof we note that the metric in $U_x \cap V$ induced from $\mathbb{R}^n$ is equivalent to the metric lifted from $B_c(0, 3r/2)$ by means of $\phi$. Therefore we can find a smaller neighborhood $N \subset N'$ of the point $x$ (depending on $\phi$) in which these two metrics are Lipschitz equivalent with the coefficient 2. Hence (3.8) will be valid for a metric ball $B \subset N$ and a measurable subset $\omega \subset B$, replacing $\phi(B)$ and $\phi(\omega)$, respectively, with $C(m) > d(m)$ replacing $d(m)$, and with the measure induced from $\mathbb{R}^n$. The proof is complete. $\square$

2. *Proof of Theorem* 1.4. According to the Bernstein doubling theorem of [FN3] and the Hadamard three circle theorem we have under the assumptions of Theorem 1.4 the following inequality.

For any compact $K \subset U$ and $F \in V_\lambda$, $\lambda \in K$,

$$(3.9) \qquad \sup_{B_c(0, r)} |F| \leq C_1 \sup_{B_c(0, 1)} |F|,$$

where $r = \frac{2+r_0}{2}$ and $C_1 = C_1(K, r_0)$. Then the function

$$F' := \frac{1}{\log C_1} \left( \log |F| - \sup_{B_c(0, r)} \log |F| \right)$$

belongs to $\mathcal{F}_r$. It remains to apply the inequality of Theorem 1.2 to obtain the required inequality

$$\sup_{B(x,\rho)} |F| \leq \left(\frac{d(n) \mid B(x,\rho) \mid}{\mid \omega \mid}\right)^{\gamma} \sup_{\omega} |F| \qquad (\rho \leq 1)$$

with $\gamma := c \log C_1$, where $c$ and $d(n)$ are the constants in Theorem 1.2. □

3. *Proof of Theorems* 1.5 *and* 1.7. Let $(Y, \mu)$ be a compact set with a positive Borel measure $\mu$. Assume that $f : Y \longrightarrow \mathbb{C}$ is a not-identically-0 continuous function satisfying

$$(3.10) \qquad \| f \|_Y := \max_Y |f| \leq \left(\frac{C\mu(Y)}{\mu(\omega)}\right)^{\alpha} \sup_{\omega} |f|$$

for every measurable set $\omega \subset Y$. Here $\alpha$ and $C$ are constants. Then $f$ satisfies

$$(3.11) \qquad \frac{1}{\mu(Y)} \int_Y \left|\log \frac{|f|}{\| f \|_Y}\right| d\mu \leq \alpha C'$$

with $C' = C'(C)$. The proof of (3.11) repeats word-for-word the arguments used in Theorem 5.1 of [Br].

We now prove Theorem 1.5. According to Theorem 1.3, for every $x$ in the $m$-dimensional compact algebraic manifold $V \subset \mathbb{R}^n$ there exists a ball $B_{r(x)}(x)$ such that

$$(3.12) \qquad \sup_B |p| \leq \left(\frac{d(m)\lambda_V(B)}{\lambda_V(\omega)}\right)^{\alpha k \deg(V)} \sup_{\omega} |p|$$

for every ball $B \subset B_{r(x)}(x)$ and polynomial $p \in \mathcal{P}_{k,n}$ with $p|_V \not\equiv 0$. Recall that $\lambda_V$ denotes the measure on $V$ induced by the Lebesgue measure of $\mathbb{R}^n$ and $B_\rho(y)$ is the metric ball of radius $\rho$ centered at $y$.

Let $\{B^l := B_{r(x_l)}(x_l)\}_{l=1}^s$ be a finite subcovering of the covering $\{B_{r(x)}(x)\}_{x \in V}$, and let $\mu$ be its multiplicity. According to the Lebesgue theorem there exists a constant $r_0 > 0$ such that every ball $B \subset V$ of radius less than $r_0$ is contained in one of the $B^l$. So inequality (3.12) holds for such a ball with the exponent $\alpha k \deg(V)$. Applying inequality (3.11) to $|p|$ and $Y := \overline{B}$ we get from (3.12)

$$(3.13) \qquad \frac{1}{\lambda_V(B)} \int_B \left|\log \frac{|p|}{\| p \|_B}\right| d\lambda_V \leq Ck \deg(V),$$

with $C = C(m, \alpha)$.

Let us now prove that (3.13) holds also for balls of radii greater than $r_0$. First, we note that for any ball $B$ of radius greater than or equal to $r_0$ there is a constant $m = m(V, r_0)$ such that

$$\lambda_V(B) \geq m.$$

Furthermore, from this it follows that

$$\frac{1}{k\lambda_V(B)} \int_B \left| \log \frac{|p|}{\|p\|_B} \right| d\lambda_V$$
$$\leq \frac{1}{km} \sum_{l=1}^{s} \left( \int_{B^l} \left| \log \frac{|p|}{\|p\|_{B^l}} \right| d\lambda_V + \log \frac{\|p\|_V}{\|p\|_{B^l}} \lambda_V(B^l) \right).$$

Applying now the Bernstein-Walsh inequality (see, e.g., [K]), we estimate the second term on the right-hand side of the inequality above as follows:

$$\frac{1}{k} \log \frac{\|p\|_V}{\|p\|_{B^l}} \leq \sup_V E_{B^l} =: C_l,$$

where $E_{B^l}$ is the global $\mathcal{L}$-extremal function of $B^l$ (see (2.26)). Since $B^l \subset V_c$ is not pluripolar, Theorem 2.2 of [S] shows that $C_l < \infty$. Putting together the two previous inequalities and (3.13) for $B^l$ we obtain

$$\frac{1}{k\lambda_V(B)} \int_B \left| \log \frac{|p|}{\|p\|_B} \right| d\lambda_V \leq \frac{\mu \operatorname{vol}(V)}{m} \left( C\deg(V) + \max_{1\leq l\leq s} C_l \right) =: C'.$$

From this inequality and (3.13) it follows that the BMO-norm of $\log|p|$ is bounded above by $C(V)k$, where $C(V) := 2\max(C', C\deg(V)) = 2C'$. The proof of Theorem 1.5 is complete. □

The proof of Theorem 1.7 is along the same lines as the previous one, so we will only give a sketch of it. Let $V$ be a compact real analytic manifold of real dimension $n$. Let $\{F_{s,\lambda}\}_{s=1}^{N}$ be a family of real analytic functions on $V$ depending real-analytically on $\lambda$ in an open subset $U \subset \mathbb{R}^m$, and let $V_\lambda := \operatorname{span}\{F_{s,\lambda}\}$. By the principle of analytic extension, for each $x \in V$ and $\lambda_0 \in K$ there exist an open neighborhood $U_x \times W_{\lambda_0,x}$ of the point $(x, \lambda_0) \in V \times U$ and an analytic embedding $\phi_{x,\lambda_0} : U_x \times W_{\lambda_0,x} \longrightarrow B_c^n(0,1) \times B_c^m(0,1) \subset \mathbb{C}^n \times \mathbb{C}^m$ such that

$$\phi_{x,\lambda_0}(U_x \times W_{\lambda_0,x}) = (B_c^n(0,1) \cap \mathbb{R}^n) \times (B_c^m(0,1) \cap \mathbb{R}^m).$$

Moreover, each $(\phi^{-1})^*(F_{s,\cdot})$, $1 \leq s \leq N$, admits a unique holomorphic extension to $B_c^n(0,1) \times B_c^m(0,1)$. Diminishing, if necessary, the neighborhood $U_x \times W_{\lambda_0,x}$ we can apply Theorem 1.4 to functions from the linear span of the extended family $\{F_{s,\cdot}^c\}$. Then returning to the coordinates on $V$ we determine that the inequality

$$(3.14) \qquad \sup_B |F| \leq \left( \frac{C(n)\lambda_V(B)}{\lambda_V(\omega)} \right)^{c_x} \sup_\omega |F|$$

holds for every $F \in V_\lambda$ with $\lambda \in W_{\lambda_0,x}$ and for every measurable subset $\omega$ of a metric ball $B \subset U_x$.

For a fixed $\lambda_0 \in K$ consider the open covering $\{U_x\}_{x \in V}$ of $V$, and select a finite subcovering $\{U_{s(\lambda_0)} := U_{x_s}\}_{s=1}^{s(\lambda_0)}$. We set now

$$\widetilde{W}_{\lambda_0} := \bigcap_{s=1}^{s(\lambda_0)} W_{\lambda_0, x_s}.$$

As in the proof of Theorem 1.5, we derive from (3.14) the inequality

$$\frac{1}{\lambda_V(B)} \int_B \left| \log \frac{|F|}{\| F \|_B} \right| d\lambda_V \leq C(V, \lambda_0)$$

for small metric balls $B$, where $F \in V_\lambda$ and $\lambda \in \widetilde{W}_{\lambda_0}$. To estimate this quantity for large balls we need only use the Bernstein doubling inequality from [FN3] instead of the Bernstein-Walsh inequality which has been used in the proof of Theorem 1.5. In this way we obtain the required estimate of the BMO-norm $|\log |F|\,|_*$ for $F \in V_\lambda, \lambda \in \widetilde{W}_{\lambda_0}$. Finally, taking a finite subcovering $\{\widetilde{W}_{\lambda_i}\}_{i=1}^s$ of the covering $\{\widetilde{W}_{\lambda_0}\}_{\lambda_0 \in K}$ and using the finite number of the estimates corresponding to $\widetilde{W}_{\lambda_i}$, $1 \leq i \leq s$, we get the global estimate of $|\log |F|\,|_*$ for any $F \in V_\lambda$ with an arbitrary $\lambda \in K$. The proof is complete. □

## 4. Concluding remarks

1. One can generalize the local inequality for polynomials of Theorem 1.3 as follows.

THEOREM 4.1. *For every regular point $x \in V$ there exists an open neighborhood $N = N(V)$ of $x$ such that for every ball $B \subset N$, measurable subset $\omega \subset B$, and polynomial $p \in \mathcal{P}_{k,n}(\mathbb{R})$ we have*

$$\sup_B |p| \leq (\alpha k \deg(V) + 1) \left( \frac{d(m) \lambda_V(B)}{\lambda_V(\omega)} \right)^{\alpha k \deg(V)} \left( \frac{1}{\lambda_V(\omega)} \int_\omega |p| \, d\lambda_V \right).$$

*Here $N$, $\alpha$, and $d(m)$ are as in Theorem* 1.3.

To prove the result it suffices to use an inequality for the distribution function of $p\,|_V$ similar to the inequality of Theorem 1.2 (see [BG] for details).

2. It was discovered in [BM] and [FN2], independently and in different ways that local Bernstein-type inequalities can be obtained from the Hadamard three circles theorem. Based on this idea, it was proved in [FN2] that Markov's inequality for algebraic functions with the constant depending on $(\deg p)^2$ linearly can be derived from the classical one-dimensional doubling Bernstein inequality. We remark that by applying a suitable version of Hadamard's theorem for the case of three polydisks one can deduce the required Markov inequality from the doubling inequality of Theorem 1.3 in a straightforward manner.

In the general case we can state that under the assumptions and notation of Theorem 1.3

$$\sup_{B} | \nabla p | \leq \frac{1}{r} \left( \frac{d(m)\lambda_V(B)}{\lambda_V(\omega)} \right)^{k\gamma} \left( \frac{1}{\lambda_V(\omega)} \int_\omega | p | \, d\lambda_V \right)$$

for a $p \in \mathcal{P}_{k,n}(\mathbb{R})$. Here $r$ is radius of $B$ and $\gamma = \gamma(N)$.

3. The BMO-properties of Theorems 1.5 and 1.7 produce a large class of examples of functions from BMO. Since a quasi-conformal change of coordinates preserves the BMO-class one can obtain additional examples.

The arguments based on the inequality of Theorem 1.2 allow us also to prove the following result:

*For any real analytic functions $f_1, \ldots, f_k$ defined on a compact real analytic manifold $V$ the function $\max_{i} \log |f_i|$ belongs to $\mathrm{BMO}(V)$.*

Using, in addition, Proposition 2.2 we can state

*For any subharmonic function $f \not\equiv -\infty$ defined in an open neighborhood of the unit circle $S^1$ the restriction $f|_{S^1}$ belongs to $\mathrm{BMO}(S^1)$.*

4. Inequality (1.2) can be used also in the problem, first posed by B. Panejah [P] for $L_2(\mathbb{R}^n)$, on the equivalence of $L_p$-norms of entire functions of exponential type over $\mathbb{R}^n$ and a relatively dense subset.

*Definition* 4.2. A measurable subset $E \subset \mathbb{R}^n$ is said to be *relatively dense* if there is a constant $L > 0$ such that

$$\delta := \inf_{x \in \mathbb{R}^n} \mathrm{mes}(B(x, L) \cap E) > 0.$$

THEOREM 4.3. *Let $E \subset \mathbb{R}^n$ be relatively dense (with constants $L$ and $\delta$) and $u : \mathbb{C}^n \longrightarrow \mathbb{R} \cup \{-\infty\}$ be a plurisubharmonic function satisfying*

$$u(z) \leq \sigma(1 + |y|) + N|x| + c$$

*for some $\sigma, N \geq 0$ and $c \in \mathbb{R}$. Then*

$$\sup_{x \in \mathbb{R}^n} u(x) \leq C'\sigma L \log \frac{c(n) L^n}{\delta} + \sup_{x \in E} u(x).$$

*Here $C'$ is an absolute constant, $z = x + iy$ and $|\cdot|$ denotes the Euclidean norm.*

Comparing this result with similar inequalities for *p-subharmonic* functions due to B. Ya. Levin and V. N. Logvinenko [LL], we note that for $n = 1$ (where both inequalities are the same) our approach leads to the improved constant $C\sigma L \log(4L/\delta)$ instead of $C\sigma L^2/\delta$ in [LL]. We have to stress, nevertheless, that our proof uses an important component of the proof in [LL] due to B. Ya. Levin.

In a forthcoming paper we present the proof of Theorem 4.3 and its generalization for regular coverings over compact algebraic manifolds (e.g. $\mathbb{R}^n$ is the covering over an $n$-torus).

The author thanks Professors Yu. Brudnyĭ and N. Levenberg for inspiring discussions.

University of Toronto, Toronto, Canada
*E-mail address*: brudnyi@math.toronto.edu